\RequirePackage{fix-cm}
\documentclass[smallextended]{svjour3}       % onecolumn (second format)
\smartqed  % flush right qed marks, e.g. at end of proof
\usepackage{graphicx}
\usepackage{booktabs}
\usepackage{longtable}
\usepackage{caption}
\usepackage{amsmath}
\usepackage{amssymb}
\usepackage{placeins}
\usepackage{subcaption}
\usepackage{listings}
\usepackage{xcolor}
\usepackage{url}

\definecolor{codegreen}{rgb}{0,0.6,0}
\definecolor{codegray}{rgb}{0.5,0.5,0.5}
\definecolor{codepurple}{rgb}{0.58,0,0.82}
\definecolor{backcolour}{rgb}{1,1,1}

\lstdefinestyle{mystyle}{
    backgroundcolor=\color{backcolour},   
    commentstyle=\color{codegreen},
    keywordstyle=\color{magenta},
    numberstyle=\tiny\color{codegray},
    stringstyle=\color{codepurple},
    basicstyle=\ttfamily\footnotesize,
    breakatwhitespace=false,         
    breaklines=true,                 
    captionpos=b,                    
    keepspaces=true,                 
    numbers=left,                    
    numbersep=5pt,                  
    showspaces=false,                
    showstringspaces=false,
    showtabs=false,                  
    tabsize=2,
    frame=single
}
\newcommand{\R}{\mathbb{R}}
\newcommand{\mb}[1]{\mathbb{#1}}
\newcommand{\mc}[1]{\mathcal{#1}}

\newcommand{\cuclarabel}{\textsc{cuclarabel}}

\newcommand{\qoco}{\textsc{qoco}}
\newcommand{\qocog}{\textsc{qoco-gpu}}

\newcommand{\gurobi}{\textsc{gurobi}}
\newcommand{\mosek}{\textsc{mosek}}

\newcommand{\cudss}{\text{cuDSS}}
\newcommand{\cuda}{\textsc{cuda}}

\usepackage[table]{xcolor}
\usepackage{calc}

\definecolor{BenchHighlight}{rgb}{0.6,0.8,0.6}
\newcommand{\winner}{\cellcolor{BenchHighlight}}
\usepackage[left=1in,right=1in]{geometry}
\usepackage{hyperref}
\hypersetup{
    colorlinks=true,
    linkcolor=blue!70,
    filecolor=PineGreen,
    urlcolor=red!50!gray!70,
    citecolor=red!50!gray!70,
    pdftitle={QOCO-GPU},
    pdfauthor={Govind M. Chari}
}

\begin{document}

\title{QOCO-GPU: A Quadratic Objective Conic Optimizer with GPU Acceleration %\thanks{This research was supported by ONR grants N000142512231 and N00014-25-1-2319.}
}
% \subtitle{Do you have a subtitle?\\ If so, write it here}

\titlerunning{QOCO-GPU}        % if too long for running head

\author{Govind M. Chari         \and
        {Beh\c{c}et}~{A\c{c}\i{}kme\c{s}e} %etc.
}


\institute{G. Chari \at
              \email{gchari@uw.edu}           %  \\
           \and
           B. {A\c{c}\i{}kme\c{s}e} \at
              \email{behcet@uw.edu}           %  \\
}

\date{Received: date / Accepted: date}

\maketitle

\begin{abstract}
We present a GPU-accelerated backend for QOCO, a C-based solver for quadratic objective second-order cone programs (SOCPs) based on a primal-dual interior point method. Our backend uses NVIDIA's cuDSS library to perform a direct sparse $LDL^\top$ factorization of the KKT system at each iteration. We also develop custom CUDA kernels for cone operations and show that parallelizing these operations is essential for achieving peak performance. Additionally, we refactor QOCO to introduce a modular backend abstraction that decouples solver logic from the underlying linear algebra implementations, allowing the existing CPU and new GPU backend to share a unified codebase. This GPU backend is accessible through a direct Python interface and through CVXPY, allowing for easy use. Numerical experiments on a range of large-scale quadratic programs and SOCPs with tens to hundreds of millions of nonzero elements in the KKT matrix, demonstrate speedups of up to 50-70 times over the CPU implementation.
% \keywords{Large-scale optimization \and GPU \and Interior point method}
\end{abstract}

\section{Introduction}
We consider the following quadratic objective second-order cone program (SOCP)

\begin{equation}\label{eq:problem}
    \begin{aligned}
        \underset{x}{\text{minimize}}
        \quad & \frac{1}{2}x^\top P x + c^\top x \\
        \text{subject to}
        \quad & Gx \preceq_\mathcal{K} h \\
        \quad & Ax = b,
    \end{aligned}        
\end{equation}
with optimization variable $x \in \R^n$. The objective is defined by $P \succeq 0$ and $c \in \R^n$. The equality and conic constraints are defined by $A \in \R^{p \times n}$, $G \in \R^{m \times n}$, $b \in \R^p$, and $h \in \R^m$. The conic inequality $Gx \preceq_\mathcal{K} h$ denotes $h - Gx \in \mathcal{K}$, where $\mathcal{K}$ is the Cartesian product 

\begin{equation*}
    \mc{K} := \mc{C}_1 \times \mc{C}_2 \times \cdots \times \mc{C}_K,
\end{equation*}
and $\mc{C}_k$ is either the non-negative orthant or a second-order cone. We assume that Problem \eqref{eq:problem} is {\bf feasible} and has a {\bf bounded} optimal objective.

In this work, we focus on solving \textit{large-scale} instances of Problem \eqref{eq:problem} with tens to hundreds of millions of nonzero elements in the KKT matrix using the {\qoco} solver \cite{chari2026qoco}, a C-based solver for quadratic objective SOCPs that implements a primal-dual interior point method \cite{vandenberghe2010cvxopt}. 

\subsection{Related work}
While GPUs have dramatically accelerated machine learning workloads, leveraging them in optimization algorithms has historically been challenging. The main difficulty is that sparse linear algebra, especially sparse matrix factorizations, is significantly more challenging to accelerate than dense linear algebra \cite{rennich2016accelerating}. As a result, most GPU accelerated optimization methods have relied on algorithms that avoid direct sparse factorizations \cite{odonoghue2021scs,schubiger2020gpu,smith2011gpu}.

Two main approaches have been considered: first-order optimization methods such as PDLP \cite{lu2025cupdlp}, which require only matrix-vector multiplications that are straightforward to accelerate on a GPU, and using indirect methods such as the preconditioned conjugate gradient method to solve linear systems rather than direct sparse factorizations \cite{odonoghue2021scs,schubiger2020gpu,smith2011gpu}. However, the first-order methods are highly sensitive to problem conditioning and can struggle to converge to high accuracy solutions. 

Interior-point methods (IPMs), however, are more robust and can achieve high accuracy, but have required using indirect methods to solve the linear system for GPU implementation \cite{smith2011gpu}. These indirect methods require only matrix-vector multiplications, which allow them to be accelerated on a GPU, but the iteration count of these methods depends on the condition number of the linear system, which can be extremely large, leading to long per-iteration runtime \cite{nocedal2006numerical}.

In late 2023, NVIDIA released {\cudss}, a high-performance direct sparse solver for GPUs, which could be integrated into optimization algorithms. Since its release, {\cudss} has been integrated as a linear system solver in the conic solvers {\cuclarabel} \cite{chen2024cuclarabel} and {\textsc{moreau}} \cite{moreau2026}, as well as in the nonlinear programming solver {\textsc{madnlp}} \cite{shin2024accelerating}.

\subsection{Contribution}
We present three main contributions. First, we develop a {\cuda} backend for the {\qoco} solver that uses {\cudss} to perform a direct sparse $LDL^\top$ factorization of the KKT system \footnote{Our implementation can be found at \url{https://github.com/qoco-org/qoco}}. This enables GPU acceleration while avoiding indirect solvers, whose performance is sensitive to conditioning. Second, we implement custom {\cuda} kernels for parallel cone operations required by the interior-point method and show that these parallelized operations are necessary to achieve good performance on the GPU. Third, we refactor {\qoco} to use a modular backend abstraction that decouples solver logic from the underlying linear algebra implementation, allowing both CPU and GPU backends to share a unified codebase.

We present benchmarks that demonstrate a $50-70 \times$ speedup over the CPU version of {\qoco} on some large problem instances. This GPU-accelerated version of {\qoco} can be called directly through a Python interface or through CVXPY \cite{diamond2016cvxpy}, making it easy to use.

Compared to {\cuclarabel} \cite{chen2024cuclarabel}, an existing GPU-accelerated IPM for conic optimization, our C-based implementation avoids the just-in-time compilation overhead of Julia, and provides a unified interface for CPU and GPU backends within Python and CVXPY scripts, making prototyping and comparing the backends easier.

\section{Implementation}
Each iteration of {\qoco}'s primal-dual IPM requires solving two linear systems with the same coefficient matrix (the KKT matrix) and performing various cone operations on $\mc{K}$. The numerical factorization of the KKT matrix can be accelerated with {\cudss} and the cone operations are parallelizable because $\mc{K}$ is the Cartesian product of many cones, $\mc{C}_k$.

\subsection{Modular backend}
To support both CPU and GPU backends without duplicating code, we refactor the main IPM loop of {\qoco} to be backend-agnostic. All solver code shared between the two backends operates on abstract data types rather than backend-specific data structures.

Specifically, we introduce the types \texttt{QOCOMatrix} and \texttt{QOCOVector}, which represent matrices and vectors used in {\qoco}. Each backend provides its own implementation of these types. The CPU implementation stores data in CPU memory, whereas the GPU implementation stores pointers to both CPU and GPU memory.

The linear system solver is abstracted through a \texttt{Linsys} interface that defines \texttt{initialize}, \texttt{factor}, \texttt{solve}, and \texttt{update} functions. Each backend provides its own implementation of this interface. The CPU backend uses the {\textsc{qdldl}} linear system solver \cite{osqp}, while the GPU backend uses {\cudss}.

Cone operations such as Jordan products, projecting onto cone $\mc{K}$, and computing the Nesterov-Todd scalings are implemented separately for the CPU and GPU backends so that these operations can be parallelized on the GPU.

At compile time, a \texttt{CMake} flag selects the desired backend and includes the corresponding implementations of \texttt{QOCOMatrix}, \texttt{QOCOVector}, \texttt{Linsys}, and the cone operations.

The Python interface exposes both backends with \texttt{pybind11}, allowing users to select the CPU or GPU implementation at runtime. This is illustrated in Listing \ref{lst:cpu-gpu-backends}.

\begin{lstlisting}[language=Python, caption=Calling CPU and GPU backends from Python interface and CVXPY, float, label=lst:cpu-gpu-backends]
# Solve with CPU backend with Python interface.
solver_cpu = qoco.QOCO(algebra="builtin")
solver_cpu.setup(n, m, p, P, c, A, b, G, h, l, nsoc, q)
result_cpu = solver_cpu.solve()

# Solve with GPU backend with Python interface.
solver_gpu = qoco.QOCO(algebra="cuda")
solver_gpu.setup(n, m, p, P, c, A, b, G, h, l, nsoc, q)
result_gpu = solver_gpu.solve()

# Solve with CPU backend in CVXPY.
problem.solve(solver="QOCO", algebra="builtin")

# Solve with GPU backend in CVXPY.
problem.solve(solver="QOCO", algebra="cuda")

\end{lstlisting}

\subsection{Linear system solver}
The primary operation accelerated by the GPU backend is the $LDL^\top$ factorization of the KKT system in Equation \eqref{eq:kkt-system}, where $W_k$ is the Nesterov-Todd scaling matrix. In the CPU version of {\qoco}, this factorization dominates the runtime and its cost grows rapidly with problem size. Therefore, we accelerate this factorization on the GPU using {\cudss}.

\begin{equation}\label{eq:kkt-system}
    \begin{bmatrix}
        P & A^\top & G^\top \\
        A & 0 & 0 \\
        G & 0 & -W_k^\top W_k
    \end{bmatrix}
    \begin{bmatrix}
        x \\
        y \\
        z
    \end{bmatrix}
    =
    \begin{bmatrix}
        r_x \\
        r_y \\
        r_z
    \end{bmatrix}
\end{equation}

{\cudss} has three phases: analysis, factor, and solve. The analysis phase is executed once and computes a fill-reducing reordering of the KKT matrix, applies this permutation to the KKT matrix, and computes a symbolic factorization that determines the sparsity pattern of the factors. This phase is the most expensive part of the solve, and for large optimization problems, it can take longer than the IPM iterations themselves. This is because {\cudss} computes the reordering on the CPU, since the graph algorithms used for this step are difficult to accelerate on the GPU. The factor phase is executed once per IPM iteration and computes the $LDL^\top$ factor of the KKT matrix. The solve phase is executed twice per IPM iteration and performs triangular solves using the factored system.

\subsection{GPU performance considerations}
Simply using {\cudss} as a linear system solver is not sufficient for developing a high-performance GPU backend due to certain implementation and hardware-level challenges that have to be addressed. Since data transfer between the CPU and GPU is significantly slower than memory access within the GPU, it is important to minimize CPU-GPU data transfers. To address this, during the setup phase of the solver, all internal {\qoco} data structures are allocated on the CPU and then copied to the GPU. During the solve phase, all computation is performed on the GPU. After convergence, the optimal solution is copied back to the CPU.

Another issue is that all cone operations on $\mc{K}$ must also be executed on the GPU, which can be extremely slow if parallelism is not exploited. As discussed earlier, operations on $\mc{K}$ decompose into independent operations on each cone $\mc{C}_k$. Therefore, these operations are trivially parallelizable. 

We have observed in our implementation that this parallelism should be exploited. In earlier versions of our GPU backend, the cone operations were executed serially on a single GPU core. Although this produced a correct algorithm, it was up to an order of magnitude slower, especially on problems with many second-order cones, and the time spent on cone operations exceeded the runtime of the matrix factorization, which should be the main computational bottleneck of the solve phase. The reason is that GPUs are designed for massively parallel workloads, while CPUs outperform GPUs on sequential computations. Therefore, parallelizing cone operations is not merely a performance optimization, but a requirement for an efficient GPU implementation. To do this, we implemented custom {\cuda} kernels for the cone operations. Each kernel implements the necessary operations for a single cone $\mc{C}_k$, and we map cones to GPU threads by launching a grid of thread blocks whose total number of threads matches the number of cones.

\section{Numerical results}

Here, we benchmark {\qocog} against the CPU implementation of {\qoco}, as well as {\cuclarabel}, {\mosek}, and {\gurobi} \footnote{Our benchmarks are publicly available at \mbox{\url{https://github.com/qoco-org/qoco-gpu-benchmarks}}}. For all solvers, we use their default settings but set the tolerances $\epsilon_{\mathrm{abs}} = \epsilon_{\mathrm{rel}} = 10^{-7}$. All results were generated on a computer with an Intel i9-14900k processor, 96 GB of RAM, and an NVIDIA GeForce RTX 5090 with 32 GB of VRAM.

In this section, problem size is defined as the total number of nonzero elements in $A$, $G$, and the upper half of $P$. We evaluate the solvers on a set of benchmark problems which include three quadratic programs (QPs): Huber regression, single-period portfolio optimization \cite{markowitz1952-ds}, and multi-period portfolio optimization \cite{boyd2017multi}, and two second-order cone programs (SOCPs): group lasso regression \cite{yuan2006grouplasso} and total variation denoising \cite{chambolle2010introduction}. The runtime of each solver, which includes setup and solve time, is limited to an hour. To compare the performance of solvers, we use performance profiles \cite{dolan2002benchmarking} and the shifted geometric mean. Details on how these metrics are computed can be found in \cite{chari2026qoco,chen2024cuclarabel,goulart2024clarabel,osqp}.

\subsection{Benchmark problems}
We consider the following QPs: single period portfolio optimization, multiperiod portfolio optimization, and Huber regression.

The single period portfolio optimization problem is 

\begin{equation*}
    \begin{split}
        \underset{x, y}{\text{minimize}} 
        \quad & x^\top D x + y^\top y - \gamma^{-1}\mu^\top x  \\
        \text{subject to} 
        \quad & y = F^\top x \\
        \quad & x \in \Delta_n \\  
        % \quad & \sum_{i=1}^{n} x_i = 1 \\  
        % \quad &  x_i \geq 0 \quad \forall i \in [1, n],\\    
    \end{split}
\end{equation*}
where $F \in \mb{R}^{k \times 100k}$ for $k$: $\{2, 5, 10, 15, 25, 50, 75, 125\}$ and $\Delta_n$ is the $n$-dimensional unit simplex.

The multiperiod portfolio optimization problem is 

\begin{equation*}
    \begin{split}
    \min_{w_t,y_t} \quad
    & \sum_{t=1}^{T}
    \left(
    w_t ^\top D w_t
    + \|y_t\|_2^2
    - \frac{1}{\gamma}\mu_t^{\top} w_t
    + \|w_t - w_{t-1}\|_2^2
    \right) \\
    \text{s.t.}\quad
    & w_0 = \bar{w}_0 \\
    & \mathbf{1}^\top w_t = 1, \qquad t=1,\dots,T \\
    & y_t = F^\top w_t, \qquad t=1,\dots,T \\
    & 0 \le y_t \le 0.01, \qquad t=1,\dots,T \\
    & \|w_t\|_1 \le L_{\max}, \qquad t=1,\dots,T,
    \end{split}
\end{equation*}
where $F \in \R^{5000 \times 50}$ is the factor loading matrix, and we solve for time horizons $T$: $\{2, 4, 6, 8, 10, 15, 20, 25, 30, 35\}$.

The Huber regression problem is given by
\begin{equation*}
    \begin{split}
        \underset{x}{\text{minimize}} 
        \quad & \sum_{i=1}^{m} \phi(a_i^\top x - b_i),  \\
    \end{split}
\end{equation*}
where $a_i^\top$ denotes the $i^{th}$ row of $A$ and $\phi: \R \to \R$ is the Huber loss. We take $A \in \mb{R}^{10 N \times N}$ for $N \in \{50, 200, 500, 1000, 2000, 4000, 6000, 10000\}$.

We consider the following SOCPs: group lasso regression and total variation denoising. 

The group lasso regression problem is
\begin{equation*}
    \begin{split}
        \underset{x}{\text{minimize}} 
        \quad & \|Ax - b\|_2^2 + \lambda \sum_{i=1}^{N} \|x^{(i)}\|_2,  \\
    \end{split}
\end{equation*}
where $x = [x^{(1)}, x^{(2)}, \ldots, x^{(N)}]$ represents a partitioning of the regression variables into groups with each $x^{(i)}$ corresponding to one group. We choose $A\in \R^{250N \times 10N}$ and solve for $N$: $\{5,20,50,100,150,300,450,750\}$.

The total variation denoising problem is 

\begin{equation*}
\begin{aligned}
\min_{U} \quad
& \mathrm{TV}(U)
+ \frac{\lambda}{2} \|U - Y\|_F^2,
\end{aligned}
\end{equation*}
where $\|\cdot\|_F$ is the Frobenius norm, $\mathrm{TV}$ is the total-variation operator and $Y$ is the corrupted image. We test with the following images from the \texttt{sk-image} collection \cite{scikit-image}: \texttt{brick}, \texttt{camera}, \texttt{grass}, \texttt{chelsea}, \texttt{coffee}, \texttt{astronaut}, \texttt{immunohistochemistry}, and \texttt{logo}.

\subsection{Benchmark results}

\begin{figure}[h]
    \captionsetup{labelfont=bf}
    \centering
    \begin{subfigure}[b]{0.49\textwidth}
        \centering
        {\includegraphics[width=\textwidth]{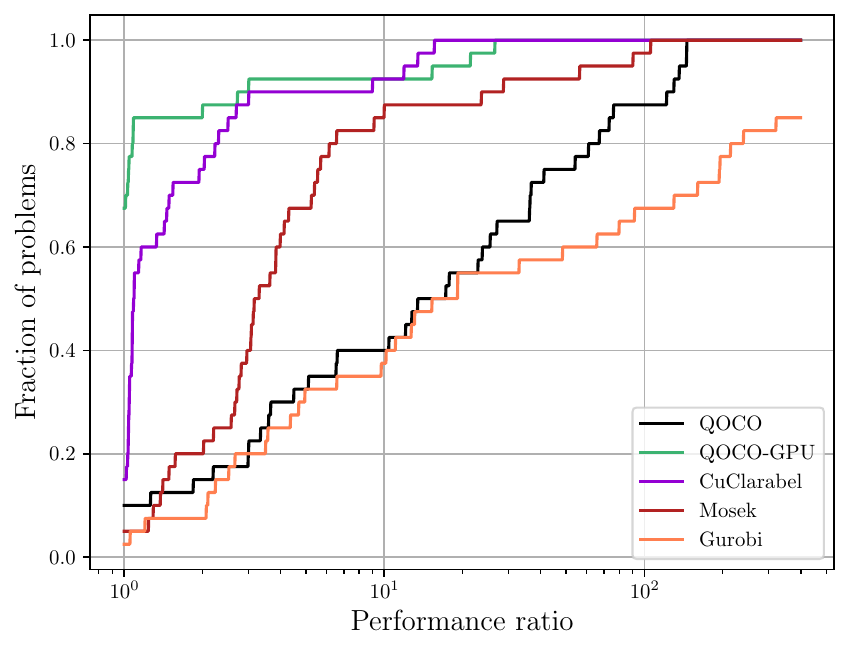}}
        \caption{Relative performance profile}
    \end{subfigure}
    \hfill
    \begin{subfigure}[b]{0.49\textwidth}
        \centering
        {\includegraphics[width=\textwidth]{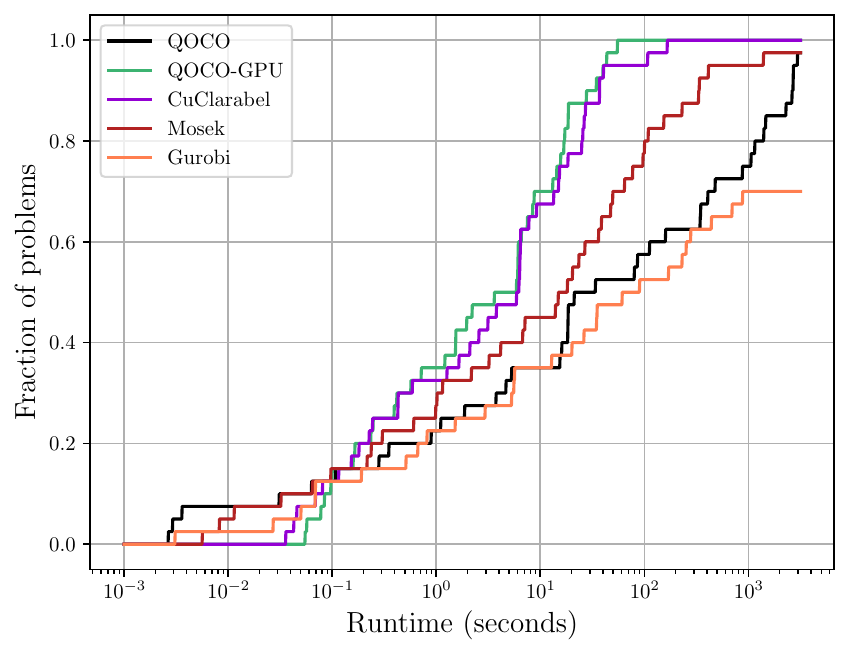}}
        \caption{Absolute performance profile}
    \end{subfigure}

    \vspace{0.5cm}

    \begin{subfigure}{1\textwidth}
        \centering
        \footnotesize
        \begin{tabular}{lccccc}
\toprule
 & QOCO-GPU & QOCO & CuClarabel & Mosek & Gurobi \\
\midrule
Shifted GM & \textbf{1.00} & 11.86 & 1.24 & 3.79 & 19.22 \\
Failure Rate (\%) & \textbf{0.0} & 2.5 & \textbf{0.0} & 2.5 & 27.5 \\
\bottomrule
\end{tabular}
        \caption{Shifted geometric means and failure rates}
      \end{subfigure}
    \caption{\bf Performance profiles for benchmark problems}
    \label{fig:benchmark-metrics}
\end{figure}

Figure \ref{fig:benchmark-metrics} shows the performance profiles and shifted geometric mean across benchmark problems and Table \ref{tab:solver_benchmarks} reports the total runtime results, including setup and solve times for each solver, where the cell corresponding to the fastest solver for each problem instance is highlighted. Table \ref{tab:solver_benchmarks} also includes percentage of the total runtime spent in setup for {\qocog} when the solver is called directly through its Python interface, as this information is not available through CVXPY. Missing entries correspond to runs that exceed the one hour time limit or did not converge. For the total variation denoising problems, {\gurobi} failed to converge to the desired accuracy. 

For smaller problems, the CPU version of {\qoco} outperforms {\qocog}, but once the KKT matrix contains around $10^5$ nonzeros, the GPU version becomes faster. For the largest problems we observe speedups of up to $70 \times$ over {\qoco}. As problems get larger, around $10^5$ to $10^6$ nonzeros in the KKT matrix, we also observe that both GPU solvers ({\qocog} and {\cuclarabel}) outperform all CPU solvers ({\qoco}, {\mosek}, and {\gurobi}). {\qocog} and {\cuclarabel} exhibit similar performance on most problems, but for a few problems, such as the largest Huber regression problem and largest group lasso regression problem, {\qocog} is approximately $2-3\times$ faster than {\cuclarabel}. 

Overall, {\qocog} has the lowest shifted geometric mean, followed by {\cuclarabel}, {\mosek}, {\qoco}, and {\gurobi}. It is expected that {\qoco} is amongst the slowest solvers on these large problems, since it uses a single-threaded matrix factorization, whereas {\mosek} and {\gurobi} use multithreaded factorizations.

Finally, for larger problems, up to $75\%$ of {\qocog}'s runtime is spent in the setup phase, where the dominant cost is the reordering step in {\cudss}'s analysis phase. This bottleneck has also been observed in \cite{pacaud2024gpu} and \cite{pacaud2024condensed}. However, the analysis phase only needs to be performed once. If the problem data changes, the reordering can be reused as long as the sparsity pattern of the problem remains fixed. This makes {\qocog} well suited for large parameteric optimization problems, where the initial cost of the analysis phase can be amortized over multiple solves.

{\footnotesize
\begin{longtable}{l r r r r r r}
\caption{\bf Runtime in seconds for benchmark problems (QOCO-GPU shows setup time percentage in parentheses)}
\label{tab:solver_benchmarks} \\

\toprule
Problem & Size & QOCO-GPU & QOCO & CuClarabel & Mosek & Gurobi \\
\midrule
\endfirsthead

\toprule
Problem & Size & QOCO-GPU & QOCO & CuClarabel & Mosek & Gurobi \\
\midrule
\endhead

\midrule
\multicolumn{7}{r}{\footnotesize Continued on next page} \\
\endfoot

\bottomrule
\endlastfoot
huber\_50 & 5500 & 0.057 (17\%) & \winner 0.003 & 0.036 & 0.011 & 0.027 \\
huber\_200 & 52000 & 0.084 (40\%) & \winner 0.031 & 0.069 & 0.097 & 0.069 \\
huber\_500 & 280000 & 0.165 (56\%) & 0.280 & \winner 0.152 & 0.605 & 0.663 \\
huber\_1000 & 1060000 & \winner 0.416 (69\%) & 1.868 & 0.431 & 4.153 & 5.278 \\
huber\_2000 & 4120000 & \winner 1.534 (78\%) & 18.510 & 1.650 & 36.219 & 20.051 \\
huber\_4000 & 16240000 & \winner 5.885 (79\%) & 159.576 & 6.435 & 331.281 & 89.658 \\
huber\_6000 & 36360000 & \winner 13.167 (78\%) & 478.210 & 25.557 & 1390.174 & 252.042 \\
huber\_10000 & 100600000 & \winner 39.816 (76\%) & 1459.514 & 107.326 & - & 439.726 \\
\midrule
portfolio\_10 & 8020 & 0.078 (11\%) & \winner 0.003 & 0.046 & 0.008 & 0.003 \\
portfolio\_50 & 140100 & 0.100 (36\%) & 0.063 & 0.115 & 0.064 & \winner 0.050 \\
portfolio\_100 & 530200 & \winner 0.159 (51\%) & 0.350 & 0.181 & 0.237 & 0.191 \\
portfolio\_200 & 2060400 & \winner 0.393 (64\%) & 5.282 & 0.427 & 1.015 & 0.813 \\
portfolio\_500 & 12651000 & \winner 2.206 (76\%) & 79.760 & 2.560 & 6.763 & 5.564 \\
portfolio\_900 & 40771800 & \winner 7.501 (71\%) & 406.031 & 7.738 & 20.330 & 26.248 \\
portfolio\_1300 & 84892600 & \winner 17.128 (67\%) & 1149.773 & 18.409 & 47.382 & 61.083 \\
portfolio\_1800 & 162543600 & \winner 34.458 (65\%) & 2618.632 & 37.001 & 108.825 & 170.539 \\
\midrule
multiperiod\_portfolio\_2 & 385600 & 0.247 & 0.891 & \winner 0.244 & 0.302 & 0.511 \\
multiperiod\_portfolio\_5 & 956500 & \winner 0.568 & 3.742 & 0.589 & 1.148 & 1.516 \\
multiperiod\_portfolio\_10 & 1908000 & \winner 1.209 & 15.418 & 1.266 & 3.215 & 5.662 \\
multiperiod\_portfolio\_15 & 2859500 & \winner 1.955 & 33.728 & 2.097 & 7.102 & 12.838 \\
multiperiod\_portfolio\_25 & 4762500 & \winner 3.601 & 85.747 & 3.776 & 14.854 & 35.032 \\
multiperiod\_portfolio\_50 & 9520000 & \winner 8.427 & 346.344 & 9.208 & 76.905 & 277.693 \\
multiperiod\_portfolio\_75 & 14277500 & \winner 14.317 & 871.908 & 15.257 & 411.356 & 693.854 \\
multiperiod\_portfolio\_125 & 23792500 & 43.295 & 2952.094 & \winner 40.338 & 229.780 & 3221.442 \\
\midrule
group\_lasso\_5 & 8805 & 0.055 (19\%) & \winner 0.004 & 0.043 & 0.006 & 0.069 \\
group\_lasso\_20 & 110220 & 0.097 (41\%) & 0.107 & 0.081 & \winner 0.032 & 2.944 \\
group\_lasso\_50 & 650550 & 0.234 (69\%) & 1.104 & 0.226 & \winner 0.216 & 34.604 \\
group\_lasso\_100 & 2551100 & \winner 0.716 (80\%) & 4.666 & 6.445 & 0.986 & 229.384 \\
group\_lasso\_150 & 5701650 & \winner 1.543 (83\%) & 16.067 & 3.134 & 2.171 & 877.229 \\
group\_lasso\_300 & 22653300 & 6.551 (82\%) & 112.120 & \winner 6.309 & 13.908 & - \\
group\_lasso\_450 & 50854950 & 15.627 (79\%) & 342.902 & \winner 14.973 & 49.441 & - \\
group\_lasso\_750 & 141008250 & \winner 54.863 (74\%) & 1400.521 & 165.127 & 336.563 & - \\
\midrule
tv\_denoising\_camera & 2092037 & \winner 6.110 & 18.302 & 6.229 & 23.394 & - \\
tv\_denoising\_grass & 2092037 & 6.064 & 21.103 & \winner 5.880 & 38.575 & - \\
tv\_denoising\_brick & 2092037 & \winner 6.141 & 18.471 & 6.381 & 18.190 & - \\
tv\_denoising\_chelsea & 2966850 & \winner 8.706 & 1058.975 & 13.401 & 26.812 & - \\
tv\_denoising\_coffee & 5267013 & \winner 16.804 & 2287.905 & 24.972 & 64.232 & - \\
tv\_denoising\_astronaut & 5753869 & \winner 18.558 & 2702.704 & 27.056 & 99.934 & - \\
tv\_denoising\_immunohistochemistry & 5753869 & \winner 18.457 & 2681.049 & 26.344 & 96.628 & - \\
tv\_denoising\_logo & 7233017 & \winner 27.720 & - & 36.955 & 153.639 & - \\
\end{longtable}
}
% \begin{figure}[h]
%     \captionsetup{labelfont=bf}
%     \centering
%     \includegraphics[width=\textwidth]{img/benchmark_problems.pdf}
%     \caption{\bf Runtime in seconds vs problem size for benchmark problems}
%     \label{fig:bench-time-vs-size}
% \end{figure}

\section{Limitations}

Although the GPU backend provides substantial speedups on large problems, there are a few limitations. The sparse direct factorization can lead to significant fill-in, which increases memory usage and can lead to out-of-memory errors on the GPU for extremely large problems. This limitation is shared with any solvers, such as {\cuclarabel}, that rely on sparse direct methods. Additionally, the GPU backend only offers a speedup for sufficiently large problems. For smaller problems, the overhead of the GPU kernel launches and the CPU-GPU memory transfers can outweigh the speedups offered by the GPU. Finally, the runtime bottleneck of this backend on large problems is the runtime of {\cudss}'s analysis phase.  

\begin{acknowledgements}
This research was supported by ONR grants N000142512231 and N00014-25-1-2319. We would like to thank Danylo Malyuta and Abhinav Kamath for their review of this paper.
\end{acknowledgements}



\FloatBarrier
\bibliographystyle{spmpsci}      % mathematics and physical sciences
\bibliography{references}   % name your BibTeX data base

\end{document}